\documentclass[12pt]{amsart}
\title[On totally real spheres in complex space]
{On totally real spheres in complex space}
\thanks{Research at MSRI is supported in part by NSF grant DMS-9022140}
\author{Xianghong Gong}
\address{Mathematical Sciences Research Institute, 
1000 Centennial Drive,
Berkeley, CA 94720}
\email{gong@msri.org}

\newtheorem{thm}{Theorem}[section]

\newtheorem{cor}[thm]{Corollary}

\newtheorem{prop}[thm]{Proposition}

\newtheorem{define}[thm]{Definition}

\newtheorem{lemma}[thm]{Lemma}

\newcommand{\cc}{{\bf C}}

\newcommand{\zz}{{\bf Z}}
\newcommand{\rr}{{\bf R}}

\newcommand{\be}[1]{\begin{equation}\label{eq:#1}}

\newcommand{\nee}{\end{equation}}

\newcommand{\ben}{\begin{equation*}}

\newcommand{\een}{\end{equation*}}

\newcommand{\nbeq}{\begin{eqnarray}}

\newcommand{\eeq}{\end{eqnarray}}

\newcommand{\bc}[1]{\begin{cor}\label{cor:#1}}

\newcommand{\ec}{\end{cor}}

\newcommand{\bt}[1]{\begin{thm}\label{thm:#1}}

\newcommand{\et}{\end{thm}}

\newcommand{\bl}[1]{\begin{lemma}\label{lemma:#1}}

\newcommand{\el}{\end{lemma}}

\newcommand{\bd}[1]{\begin{define}\label{lemma:#1}}

\newcommand{\ed}{\end{define}}

\newcommand{\bp}[1]{\begin{prop}\label{prop:#1}}

\newcommand{\ep}{\end{prop}}

\newcommand{\jq}[1]{\mathopen{<} #1\mathclose{>}}%produce <*,*>

\newcommand{\rl}[1]{Lemma~\ref{lemma:#1}}

\newcommand{\nrc}[1]{Corollary~\ref{cor:#1}}

\newcommand{\re}[1]{(\ref{eq:#1})}

\newcommand{\rp}[1]{Proposition~\ref{prop:#1}}

\newcommand{\rt}[1]{Theorem~\ref{thm:#1}}

\begin{document}

\begin{abstract}
We shall prove that there are totally real and real analytic embeddings of 
$S^k$ in $\cc^n$ which are not biholomorphically equivalent if $k\geq 5$ and
$n=k+2[\frac{k-1}{4}]$. We also show that
a smooth manifold $M$ admits a totally real immersion in $\cc^n$ with
a trivial complex normal bundle 
if and only if the complexified tangent bundle
of $M$ is trivial. The latter is proved by applying Gromov's weak homotopy 
equivalence principle for totally real immersions to Hirsch's
transversal fields theory.
\end{abstract}

\maketitle

\setcounter{thm}{0}\setcounter{equation}{0}

\section{Introduction}

Let $M, N$ be two totally real and real analytic submanifolds in $\cc^n$.
 We say that $M$ and $N$ are {\it biholomorphically equivalent} if
  there is a biholomorphic mapping $F$ defined in a neighborhood of $M$
 such that $F(M)=N$.
As a standard fact of
complexification, one knows that all totally real and real analytic
embeddings of $M$ in $C^n$ are biholomorphically equivalent if $M$ is of
maximal dimension $n$. However, the topology of  the manifold plays a major
role in the existence of totally real immersions or embeddings. For instance,
R.O.~Wells~\cite{wells} proved that if an $n$-dimensional
compact and orientable manifold $M$ admits
a totally real
embedding in $\cc^n$, then its Euler number must vanish.
It was also observed by Wells that if  
 $M$ is a manifold of  dimension  $n$ and it admits a totally
real immersion in $\cc^n$, then its complexified tangent
bundle $T^cM=TM\otimes\cc$ is trivial.
Conversely,  the triviality of $T^cM$ also implies the existence of 
totally real immersions of $M$ in $\cc^n$.
This was obtained by M.L.~Gromov in~\cite{gromov} through the method of convex
integration.

The sphere $S^
k\colon x_1^2+\ldots+x_{k+1}^2=1 $ in $\rr^{k+1}$ gives us  
 a trivial totally real embedding of
$S^k$ in
$\cc^{k+1}$. On the other hand, the works of Gromov~\cite{gromov},
 Ahern-Rudin~\cite{ahern} and Stout-Zame~\cite{stout}
tell us that $S^k$ admits a totally real and real analytic embedding in  $\cc^k$ if and only if $k=1, 3$.
Our main result is the following.
\bt{1}
If $k\leq 4$, all totally real and real analytic embeddings of $S^k$ in
$\cc^n$ are biholomorphic equivalent.
If $k\geq 5$ and
$n_k=k+2[\frac{k-1}{4}]$, there exist totally real and real
analytic
embeddings of $S^k$ in $\cc^{n_k}$ which are not biholomorphically equivalent,
 while all totally  real and real analytic embeddings of $S^k$ in
 $\cc^n$ are biholomorphically equivalent if $n>n_k$.
 \et
In fact, we shall prove a slightly stronger result that all totally real and
real analytic embeddings of $S^k$ in $\cc^n$ are unimodularly equivalent
 if $k\leq 4$ and $n>k$, i.e. they are biholomorphic equivalent through 
a  
 mapping $F$   
which preserves the holomorphic $n$-form $dz_1\wedge\ldots\wedge dz_n$.
 We should mention that using a transversality argument~\cite{fr},
F.~Forstneri\v{c} and J.-P. Rosay
showed that
for a  compact manifold $M$  of dimension $k$,
all its
 totally real and real
analytic embeddings in $\cc^n$  are biholomorphically equivalent if
$n\geq 3k/2$.

The proof of \rt{1} is not constructive. It depends on the
weak homotopy equivalence (w.h.e.) principle for totally real immersions established 
by Gromov in~\cite{gromov}.
To show \rt{1}, we also need to understand the role which the normal
bundle of totally real immersions plays.
Recall that a $C^1$-smooth mapping
$f\colon M\to \cc^n$ is a {\it totally real immersion} if 
$f_*T_xM$
spans a $k$-dimensional complex linear subspace of $T_{f(x)}\cc^n$ for each $x\in M$.
We define the {\it complex} normal bundle
of the immersion $f$, denoted by $\nu_f$, to be the complex vector bundle
whose fiber over $x\in M$ is the quotient of $T_{f(x)} \cc^n$ by the
complex linear span of $f_*T_xM$.
We shall see that two totally real and
real analytic embeddings of a manifold are biholomorphically equivalent
if and only if their complex normal bundles are topologically equivalent.
The normal bundle of immersions plays quite important role in the works
of S.~Smale~\cite{smale} and M.W.~Hirsch~\cite{hirsch}.
Analogous to the results of Hirsch~\cite{hirsch} about transversal fields of
smooth immersions, we obtain the following.
\bt{2}
Let $f\colon M\to\cc^n$ be a $C^1$ totally real immersion. Assume that the
complex normal bundle $\nu_f$ has a topologically trivial subbundle of
rank $r$. Then
there is a regular homotopy $f_t \colon M\to\cc^n$ of $C^1$ totally real
immersions such that $f_0=f$ and $f_1\colon M\to \cc^{n-r}$.
 \et

We now draw some conclusions from \rt{2}.
\bc{1}
Let $M$ be a smooth  manifold of dimension $k$. Then $M$ admits a totally real immersion
in $\cc^n$ with a trivial complex normal bundle if and only if
there exists a totally real immersion
of $M$ in $\cc^k$, i.e.
the complexified
tangent bundle $T^cM$ is trivial.
\ec
A smooth manifold $M$ is said to be {\it stably parallelizable} if the 
tangent bundle of 
$M\times \rr$ is trivial.
For instance,
the boundary of a smooth domain
in euclidean space is always stably parallelizable.
By a theorem of Hirsch~\cite{hirsch},
$M$ is stably parallelizable if and only if it is orientable and admits
an immersion
in $\rr^{n+1}$.
>From \rt{2}, we have
the following.
\bc{2}
Let $M$ be a  manifold of dimension $n$ which is immersible in $\rr^{n+1}$.
Then
$T^cM$ is trivial if $M$ is orientable, or $M$ is non-orientable
with
$H^2(M,\zz)=0$.
 \ec
In fact, Gromov  proved a stronger result that $M$ admits an exact Lagrangian
immersion in euclidean space when $M$ is stably parallelizable (see~\cite{gromovbk}, p.~61).
%check forst's paper, may state for surface, any case the next remark need
%to restate
%
%It's proved by F that orientable surface has totall real immersion in C^2
% For the non-orientable surface M, H^2=\zz_2. cor is incolusive. In
% stead, from F's result and FR's result we know that the normal line bundle
% is trivial if and only if its genus is odd.
%
One notices that all real surface $M$ can be immersed in
$\rr^3\subset\cc^3$.
On the other hand,  Forstneri\v{c}~\cite{forstneric}  proved that a
non-orientable
compact surface admits a totally real immersion in $\cc^2$  if and only if
its genus is even. Therefore,
the condition that $H^2(M,\zz)$ vanishes is essential in \nrc{2}.

The paper is organized as follows. In section two we shall discuss
Gromov's w.h.e.-principal for totally real immersions. Section three is
devoted to the proof of \rt{2} and \nrc{2}. The proof of \rt{1} will be
given in the last section, where we shall also make essential uses of
homotopy groups of complex Stiefel manifolds  obtained by
M.L.~Kervaire~\cite{kervaire} and F.~Sigrist~\cite{sigrist}.

\setcounter{thm}{0}\setcounter{equation}{0}
\section{Classification of totally  real immersions}
In this section, we shall first recall Gromov's w.h.e.-principal for ample
differential relations established in~\cite{gromov}. We shall also
discuss the group structure on the regular homotopy classes of totally
real immersions of $S^k$ in $\cc^n$.

Let $M$ be a smooth manifold of dimension $k$.  Assume that $k\leq n$.
  By a {\it regular homotopy} $f_t $ of totally real
$C^1$-immersions
of $M$ in $\cc^n$, one means that for each $t\in [0,1]$, $f_t$ is a 
totally real $C^1$-immersion,
and $df_t\colon TM\to T\cc^n$ depends on $t\in [0,1]$ continuously.
Mappings from $M$ to $\cc^n$ can  be identified
with sections of the
 trivial bundle $X=M\times\cc^n\to M$.
If $f\colon M\to \cc^n$ is a $C^1$-mapping defined in a neighborhood
of $x\in M$, we define the 1-jet 
 of $f$ at $x$ to be $J_x^1f=(f(x), df_x)$.
We denote by $X^1$
the space of 1-jets of $C^1$ sections of $X\to M
$.
Then it is easy to see that $X^1$ is a
fibration over $M$ whose  fiber
over $x\in M$ consists of
$\rr$-linear mappings from $T_xM$ to $T_z\cc^n$ for some $z\in\cc^n$.
We shall adapt (compact-open) $C^0$-topology on $X^1$. Here, we should say
a few words about the topologies used in the sequel. For a homotopy of sections
or immersions, we shall always use the (compact-open) $C^r$-topology, i.e.
the weak $C^r$ topology. For approximating a mapping or function, we shall
always use the fine $C^r$-topology. The reader is referred to~\cite{hirschbk}
for basic properties of these two kinds of topologies. 

By ${\scriptsize\Sigma}_x$, one denotes the set of 1-jets $J_x^1f$ such that
$df_x(T_xM)$ spans a complex linear subspace of $T_{f(x)}\cc^n$ of rank
$k$. This is equivalent to say that 
the complexification $df_x\otimes \cc\colon T_x^cM\to
T_{f(x)}\cc^n$ is injective.
Let ${\scriptsize\Sigma}$ be the union of
$\Sigma_x$ for all $x\in M$.  
Then $\Omega=
X^1\setminus { \Sigma
}$ is an open subset of $X^1$, which is
called a {\it
totally real differential relation}.
Thus, a $C^1$-mapping $f\colon M\to\cc^n$
is totally real if and only if $J^1f$ maps $M$ into
$\Omega$.

For each $x_0\in M$, we choose local coordinates $u_1,
 \ldots, u_k$
in a neighborhood $U$ of $x_0$. Fix $x\in U$ and $1\leq j\leq n$.
Let $Z$ be the set of 1-jets $J_x^1f$ satisfying
$$
f(x)=z,\quad f_{u_i}(x)=v_i,\quad i\neq j,
$$
where $z$ and $v_i (i\neq j)$ are fixed $k$ vectors in $\cc^n$. Then either
$Z
\setminus
{ \Sigma}$ is  an empty set, or the linear convex hull of
each  connected component of
$Z\setminus{ \Sigma}$ is the whole affine space $Z$. According to
the terminology of Gromov~\cite{gromov}, $\Omega$ is said to be {\it ample}
in the coordinate directions.
To see this, we notice that
if $v_1,\ldots, \hat v_j, \ldots, v_k$ are not $\cc$-linearly independent,
then $Z\subset{ \Sigma}$.
Otherwise, $Z\cap{ \Sigma}$ consists of all 1-jets  $(z, v_1,
\ldots, v, \ldots, v_k)$ such
that $v$ is a $\cc$-linear combination of
$v_1,\ldots, \hat v_j, \ldots, v_k$. Therefore, $Z\cap  {\Sigma}$ is a subspace of the affine space $Z$ with real codimension
$2n-2(k-1)\geq 2$. This implies that $Z\setminus{ \Sigma}$ is connected and it spans
the whole space $Z$.
As a consequence of Theorem~$1.3.1$ in~\cite{gromov}, we can
state the following result.
\begin{thm}[\text{Gromov},~\cite{gromov}] \label{thm:g}
Let $M$ be a smooth manifold of dimension $k\leq n$, and let $\Omega$ and $J^1$
be 
as above.
Then $J^1$ is a one-to-one correspondence between the regular homotopy classes
of totally real $C^1$-immersion of $M$ in $\cc^n$ and the homotopy classes of
continuous sections of $\Omega\to M$.
 In particular, $M$ admits a totally real
immersion in $\cc^n$ if and only if  $\Omega\to M$ has
a global continuous section.
\et

We now consider the case that 
$M$ is the sphere $S^k$. Here we need the fact that
the complexified tangent bundle $T^cS^k$ of $S^k$ is trivial.
This follows from the existence of totally real immersion of $S^k$
in $\cc^k$. An explicit example of Lagrangian (whence totally real) immersions
 of $S^k$ was
constructed
by A.~Weinstein~\cite{weinsteinbk}.  As we mentioned earlier, Gromov showed that
all stably parallelizable  manifold, such as a sphere,  admits an exact Lagrangian
immersion (see~\cite{gromovbk}, p.~61). It seems to us that there are no other proof in the
literature about
the triviality of $T^cS^k$.
Throughout the whole paper,
 we shall fix a topological trivialization of
$T^cS^k=S^k\times\cc^k$.

Recall that the complex Stiefel manifold $V_{n,k}$ consists of
$k$-frames of $\cc^n$, i.e. the space of ordered $k$ 
linear independent vectors in $\cc^n$. With the fixed trivialization for $T^cS^k$,  the global sections
of $\Omega\to M$ can be identified with mappings from
$S^k$ to $\cc^n\times V_{n,k}$ as follows. Let $e_1,\ldots, e_
k$ be
the set of $\cc$-linearly independent continuous sections
which defines the trivialization of $T^cS^k$.
 Assume that $\phi\colon M\to \Omega$ is a global
section. Then $\phi=(f,\varphi)$ where $f\colon M\to \cc^n$ and
$\varphi\colon TM\to T\cc^n$ satisfy the property that for each $x\in M$,
$\varphi(x)\colon T_xM \to T_{f(x)}\cc^n$ is $\rr$-linear
and its complexification  is injective. Hence,
$$
v(x)=(\varphi(x)(e_1(x)),\ldots, \varphi(x)(e_k(x)))
$$ 
is a set of linearly
independent $k$ vectors of $T_{f(x)}\cc^n$.
Denote by $U_{n,k}$  the space of unitary $k$-frames of $\cc^n$, where by
a unitary $k$-frame $(v_1,\ldots, v_k)$, one means that $v_1,\ldots,v_k$ satisfy the
condition $\jq{v_i,v_j}=\delta_{i,j}$ for the standard hermitian inner
product $\jq{\cdot,\cdot}$ on $\cc^n$.
>From the well-known normalization, one knows
 that $V_{n,m}$ is the product of $U_{n,m}$ with  the space $T$ of
 upper-triangle matrices with positive eigenvalues. Thus, 
 a homotopy $\phi_t\colon S^k\to\Omega$ induces a homotopy 
$$
(f_t, A_t, \varphi_t)\colon S^k\to \cc^n\times T\times U_{n,k}. 
$$
 Since
 $\cc^n$ and $T$ are contractible, we see that the set 
of homotopy classes of
 sections
 of $
\Omega\to S^k$ is the same as the homotopy classes of  mappings
 from $S^k$ to $U_{n,k}$. It is well-known that the homotopy classes of 
mappings
 from $S^k$ to $U_{n,k}$ is the homotopy group $\pi_k(U_{n,k})$
 (see p.~88 in~\cite{steenrod}, p.~211 in~\cite{bottbk}).
Thus, we obtain a one-to-one
 mapping
$$
  j_*\colon I(S^k,\cc^n)\to \pi_k(U_{n,k}),
$$
where $I(S^k,\cc^n)$ stands for the set of regular homotopy classes of
totally real immersions of
 $S^k$ in $\cc^n$.

 With the fixed trivialization of $T^cS^k$,
 $j_*$  is a {\it canonical} mapping in the sense that there is
a commutative diagram
\be{can} \begin{array}{cccc} %line one
I(S^k,\cc^n) &  \stackrel{j_*}{\longrightarrow}   &  \pi_k(U_{n,k})
 & \\ %line two
  {\downarrow} i_*  & &
 {\downarrow} i_*' & \\ %line three
I(S^k,\cc^N) &  \stackrel{j_*'}{
\longrightarrow}   &  \pi_k(U_{N,k}),&  N>n,
\end{array} \end{equation}
 where the homomorphism $i_*$ is induced by the inclusion of
sending a totally real immersion
$f\colon S^k\to\cc^n$ to a totally real immersion    $(f,0)\colon
S^k\to\cc^N$,
and $i_*'$ is induced by regarding a unitary $k$-frame of $\cc^n$ as a
unitary $k$-frame of $\cc^N$.
In other words,
the adapted group structure on $I(S^k,\cc^m)$ is preserved under the
inclusion $\cc^n\subset\cc^N$. This will be important for us to prove
\rt{1}.

\setcounter{thm}{0}\setcounter{equation}{0}
\section{Proof of \rt{2}}
In this section we shall apply Gromov's w.h.e.-principle  to prove
\rt{2}.
With necessary modifications, we shall follow
 very closely  the
proof of Theorem~6.4 in~\cite{hirsch}.

Let $M$ be a  manifold of dimension $k$.
By $FM$, one denotes the (complexified) $k$-frame bundle of $M$, which consists of
 ordered $k$ linearly independent
vectors
in $T^cM$ with the same base point.
For each $e=(v_1,
\ldots,v_k)\in F_xM$ and $g=(g_{i,j})\in GL(k,\cc)$, we define $g\cdot e$
to be the $k$-frame $(v_1^\prime, \ldots, v_k^\prime)$ with
$v_i^\prime=\sum_{j=1}^kg_{i,j}v_j$.
Then $FM$ is a principal $GL(k, \cc)$-bundle
over $M$.
 For $k< m\leq n$, we denote by
$E_{n,m}$ be the
associated bundle of $FM$ with fiber $V_{n,m}$.
More precisely, we define a $GL(k,\cc)$-action on $V_{n,m}$ by 
$$
(g,(v_1,\ldots,v_m))\mapsto (g\cdot(v_1, \ldots, v_k), v_{k+1},\ldots,
v_m).
$$
 Then  $E_{n,m}$ is the set of equivalence classes  of the relation
$\sim$ on $FM\times V_{n,m}$ with $(e,f)\sim (g\cdot e, g\cdot f)$ for
all $g\in GL(k,\cc)$.
The bundle projection $p_{n,m}\colon E_{n,m}\to M$ is induced by
the composed projection $FM\times V_{n,m}\to FM\to M$.
Let $p_{n,k}^{n,m}\colon
V_{n,m}
\to V_{n,k}$ be the projection of deleting the last $m-k$ vectors
from each $m$-frame.
Notice that the $GL(k,\cc)$-action on $V_{n,m}$ does not affect the last
$m-k$ vectors of an $m$-frame of $\cc^n$. Hence,
  $p_{n,k}^{n,m}$ induces a projection from $E_{n,m}$ to
$E_{n,k}$ such that
$p_{n,m}=p_{n,k}\circ p_{n,k}^{n,m}$.
 We further remark that
$p_{n,k}^{n,m}\colon
E_{n,m}\to E_{n,k}$
is a fiber bundle. 
 In particular, $p_{n,k}^{n,m}\colon
E_{n,m}\to E_{n,k}$ has  the covering  homotopy property, i.e. for any
finite polyhedron $P$, a homotopy
$h_t\colon P\to E_{n,k}$ has a lifting $\tilde h_t\colon P\to E_{n,m}$, if
the initial lifting $\tilde h_0$ exists.

To use the covering homotopy property, we consider the set of sections of
the fiber
bundle $E_{n,m}\to M$.  Given a section $s\colon M\to E_{n,m}$, we define a
mapping $\varphi\colon FM\to V_{n,m}$ by  $\varphi(e)=f$ if $s(x)$ is the
equivalence class of $(e,f)\in FM\times V_{n,m}$.  Then $\varphi$ is
well-defined.
For if $(e,f)$ and $(e,f^\prime)$ are in the same equivalence
 class $s(x)$.
Then there is $g\in GL(k,\cc)$ such that $g\cdot e=e$ and $g\cdot f=
f^\prime$.
 Obviously, $g\cdot e=e$ implies that
$g=\mbox{id}$. Hence, $f^\prime=f$. Moreover, if $(e,f)$ is in the
equivalence
class $s(x)$, so is $(g\cdot e, g\cdot
f)$. Hence, $\varphi(g\cdot e)=
g\cdot\varphi(e)$, i.e. $\varphi$ is a $GL(k,\cc)$-{\it equivariant} mapping.
Conversely, given a $GL(k,\cc)$-equivariant mapping $
\varphi\colon FM\to V_{n,m}$, we set $s(x)$ to be the equivalence class of
$(e,\varphi(e))$
 in $E_{n,m}$ for $e\in F_xM$. Thus, there is a one-to-one
correspondence between the set of $GL(k,\cc)$-equivariant mappings from $FM$
to $V_{n,m}$ and the set of sections of the fiber bundle $E_{n,m}\to M$.
Therefore, the covering homotopy property
gives us the following.
\bt{lift} %cite hirsch
Let $\varphi_t$
be a homotopy of
$G(k,\cc)$-equivariant mappings from $FM$ to $V_{n,k}$. Assume that $\varphi_0$ has a lifting
$\tilde\varphi_0\colon FM\to V_{n,k+r}$
 such that $\tilde\varphi_0$ is $GL(k,\cc)$-equivariant. Then there is
a homotopy $\tilde \varphi_t$ of lifting of $\varphi_t$ such that
each   $\tilde\varphi_t\colon FM\to V_{n,k+r}$  is
 $GL(k,\cc)$-equivariant.
\et
We now let
 $p'\colon V_{n,k+r}\to V_{n,r}$
 be the mapping of projecting a $(k+r)$-frame
 to its last $r$ components.
 Let   $\phi\colon FM\to V_{n,k+r}$ be a $GL(k,\cc)$-equivariant mapping.
 Then $\varphi=p\phi\colon FM\to V_{n,k}$ 
  is also a $GL(k,\cc)$-equivariant mapping.
  Notice that $GL(k,\cc)$ acts
 on each fiber of $FM$ transitively and that the last $r$ components of a
 $(k+r)$-frame of $\cc^n$ is fixed under the $GL(k,\cc)$-action. This implies
 that
$p'\varphi$ is constant along fibers of $FM$. Therefore, $p'\varphi$
 is the lifting of some mapping $\psi\colon M\to V_{n,r}$.
 Following~\cite{hirsch}, we shall call $\psi\colon M\to V_{n,r}$ a {\it
 transversal}
 $r$-field of $\varphi\colon M\to V_{n,k}$. 
Thus, we identify  the set of $GL(k,\cc)$-equivariant
 mappings from $FM$ to $V_{n,k+r}$
 with the set of $GL(k,\cc)$-equivalent mappings from $FM$ to
$V_{n,k}$
with  transversal $r$-fields.
In general, we say that
$\psi\colon M\to V_{n,r}$ is {\it transverse} to a totally real immersion $f\colon M\to \cc^n$ if
for each $x\in M$ and
 $\psi(x)=(v_1,\ldots, v_r)$,  $v_j$ is not contained in 
$f_*(T_x^cM)\ (1\leq j\leq r)$. 
\bc{normal}
Let $f_t$ be a homotopy of totally real immersions of $M$ in $\cc^n$. 
Assume that
$\nu_{f_0}$ has a topological trivial subbundle
of rank $r$. 
Then $\nu_{f_1}$ also contains a topological trivial subbundle
of rank $r$. 
\ec
\begin{proof} We identify $f_{t*}\colon T^cM\to T\cc^n$ with a homotopy $f_t$ of mappings from $M$ to $\cc^n$ and a homotopy $\varphi_t$ of $GL(k, \cc)$-equivariant mappings from $FM$ to $V_{n,k}$. By assumptions, $\varphi_0$ has a transversal $r$-field. Hence, \rt{lift} implies that $\varphi_1$ also admits a transversal $r$-filed, i.e. $\nu_{f_1}$ has a trivial subbundle of rank $r$. 
\end{proof}
We also need the following.
\bl{smooth}
Let $M$ be a smooth manifold and $f\colon M\to\cc^n$
 a $C^\infty$-smooth totally real immersion. Assume that 
$\nu_f$ has  a topologically trivial
subbundle
 of rank $r$. Then there is a smooth unitary $r$-field
$\psi\colon M\to U_{n,r}$ which is transverse to $f$. 
\el
\begin{proof} By assumptions, there is a continuous $r$-field $\psi_0\colon M\to 
V_{n,r}$
which is transverse to the immersion $f$. Using the approximation 
in fine $C^0$-topology, 
we can replace $\psi_0$ by a smooth $r$-field
$\psi_1$ which is still transverse to the immersion $f$. 
We now project
$\psi_1(x)$
to the orthogonal complement of $f_*T_x^cM$. By the well-known
normalization,  
we readily obtain the desired unitary $r$-filed.
\end{proof}

%\vspace{1ex}
\noindent {\bf Proof of \rt{2}}. %following hirsch
 Let $f\colon M\to\cc^n$ be a $C^1$-smooth totally real immersion. Assume that the complex
 normal bundle $\nu_f$ has a topologically trivial subbundle of rank $r$. 
We shall seek a
 homotopy $f_t$ of totally real immersions of $M$ in $\cc^n$ such
 that $f_0=f$ and $f_1(M)\subset \cc^{n-1}$. We further require that the
 complex normal bundle of $f_1\colon M\to\cc^{n-1}$
has a trivial subbundle
 of rank $r-1$. Thus, \rt{2} follows from the 
 induction. 
  
We first find a smooth mapping $g\colon M\to\cc^n$ which is sufficiently
close to $f$ in
fine $C^1$-topology such that each
 $f_t=(1-t)f+tg$ is still a totally real immersion for $0\leq t\leq 1$.
Rename $g$ by $f$. Then, \nrc{normal} implies that $\nu_f$ still contains
a trivial subbundle of rank $r$. By \rl{smooth}, there exists
 a smooth unitary
$r$-field
$\psi=(\xi_1,\ldots, \xi_r)$ which is transverse to
the smooth totally real immersion  $f\colon M\to \cc^n$.

>From the standard embedding $S^
{2n-1}\subset\cc^n$, one gets a
complex vector bundle $T^{(1,0)}S^{2n-1}$ whose fiber
over $v\in S^{2n-1}$ consists of vectors in $\cc^n$
which are orthogonal to $v$ with respect to the standard hermitian metric
on $\cc^n$. Let $E$ be the frame
bundle over $S^{2n-1}$ whose fiber
 consists of linear independent vectors $v_1,\ldots, v_{k+r-1}$
  of $T^{(1,0)}S^{2n-1}$
 with the same base point.
Consider the  mapping
 \be{f}
 \tilde \xi_r(x)\colon e\mapsto (df_x(e), \xi_1(x),\ldots, \xi_{r-1}(x)),
 \quad e\in F_xM.
 \end{equation}
Since $\xi_r$ is orthogonal to $f_*T^cM$ and $\xi_j\ (j<r)$, then
$\tilde\xi_r(x)$ maps $F_xM$ into $E_{\xi_r(x)}$. Hence, $\tilde\xi_r\colon
FM\to E$ is a $GL(k,\cc)$-equivariant mapping which
 covers the mapping $\xi_r\colon M\to S^{2n-1}$.
Since the dimension
of $M$ is less than $2n-1$, the smoothness of the mapping
$\xi_r$ implies that
there is $y_0\notin  \xi_r(M)$. Set $Y=S^{2n-1}
-\{y_0\}$. Since $Y$ is contractible, then there
is a homotopy $h_t\colon M\to Y$ such that $h_0=\xi_r$ and
$h_1\equiv y_1
\in Y$. Also, there is
a trivialization
 \be{trivial}
  T^{(1,0)}S^{2n-1}|_Y=Y
\times \tilde\cc^{n-1}, \quad \tilde\cc^{n-1}=T_{y_1}^{(1,0)}S^{2n-1}.
\end{equation}
Therefore, $\tilde\xi_r$ can be written as $(h_0, \phi_0)$ with 
$\phi_0$ a
$GL(k,\cc)$-equivariant mapping from $FM$ to $E_{y_1}$.
Put $\tilde h_t=(h_t,\phi_0)$.
Returning to $T^{(1,0)}S^{2n-1}|_Y$ from the trivialization \re{trivial}, 
we obtain a homotopy $\tilde h_t$  of $GL(k,\cc)$-equivariant mappings
 from $FM$ to $E|_Y$. 
Returning
to the ambient space $\cc^n$, we then have
$T^{(1,0)}S^{2n-1}\subset S^{2n-1}\times\cc^n$. Thus, we obtain
 a homotopy $\tilde h_t$ of $GL(k,\cc)$-equivariant mappings from $FM$
to $\cc^n\times V_{n,k+r-1}$ such that $\tilde h_1\colon FM\to E_{y_1}$. 

Let us identify $y_1$ with a point 
in $\cc^{n-1}$. Put $T_{y_1}\cc^{n-1}=\tilde\cc^{n-1}$. Then we have 
$E_{y_1}=y_1\times V_{n-1,k+r-1}$. 
We now write $\tilde h_t$ as $(h_t, \varphi_t, \psi_t)$, where each $\varphi_t\colon 
FM\to V_{n,k}$  is $GL(k,\cc)$-equivariant, and  $\psi_t$ is a transversal 
$(r-1)$-field of $\varphi_t$. From \re{f}, it is clear that $\varphi_0
\colon FM\to V_{n,k}$
is just the $GL(k,\cc)$-equivariant mapping induced by 
$f_*\colon T^cM\to T\cc^n$. This implies that $f_*$ and $(h_1,\varphi_1)$
are homotopic as fiberwisely injective $\cc$-linear mappings from $T^cM$
to $T\cc^n$. 
Notice that $h_1\equiv y_1\in\cc^{n-1}$ and 
$\varphi_1\colon FM\to V_{n-1,k}$.
Now \rt{g} implies that there is a totally real
immersion $g\colon M\to \cc^{n-1}$
such that
$g_*$ and $(h_1,\varphi_1)$ are joined by a homotopy $k_t$ of 
fiberwisely injective $\cc$-linear mappings from $T^cM$ to $T\cc^{n-1}$.
Thus, we have proved that $f_*\colon T^cM\to T\cc^n$
is homotopic to $g_*\colon T^cM\to T\cc^{n-1}$ as fiberwisely 
injective $\cc$-linear
mappings from
$T^cM$ to $T\cc^n$. Using \rt{g} again, we know that there is a homotopy of
totally real immersions joining $f$ and $g$.
To complete the proof of \rt{2}, we need to show that
the complex normal bundle of $g\colon M\to \cc^{n-1}$
has a topologically trivial subbundle of rank $r-1$.
To this end, we notice that $\varphi_1$ has a transversal $(r-1)$-field 
$\psi_1$ in $\cc^{n-1}$
since $\tilde h_1
\colon FM\to E_{y_1}\equiv V_{n-1,k+r-1}$.
By applying \nrc{normal} to the homotopy $k_t$, we see that the complex
normal bundle of $g$ has a trivial subbundle of rank $r-1$. 
The proof of \rt{2} is complete.

\vspace{1ex}
\noindent {\bf Proof of \nrc{2}.}
Assume that $M$ is
immersed in $\rr^{n+1}\subset\cc^{n+1}$. When
$M$ is orientable, the real normal line bundle of the immersion of $M$
in $\rr^{n+1}$ is trivial.
Let $\nu$ be a normal
vector field of $M$. It is obvious that $\nu$ is still transverse to
$T^cM$. Therefore, \rt{2} implies that $M$  admits a totally real
immersion in
$\cc^n$,
i.e. $T^cM$ is trivial. We now consider the case that $M$ is non-orientable
and $H^2(M,\zz)=0$. By a theorem of Whitney, a smooth manifold has a compatible
analytic structure. Using the analytic approximations of smooth mappings in fine
$C^r$-topology~\cite{grauertlevi},  we can replace the original immersion of $M$ in $\rr^{n+1}$ 
by an analytic immersion $f\colon M\to\rr^{n+1}$. Now the complexification
of $f$ gives us a holomorphic mapping from $M^c$ to $\cc^{n+1}$, where $M^c$
is a complexification of $M$ depending on $f$. 
By Grauert's Tube theorem~\cite{grauertlevi}, there is a Stein neighborhood $X$ of $M$ 
in $M^c$ such that $M$ is a strong deformation retraction of $X$.
Since $X$ is Stein,
then all
holomorphic line bundle over $X$ is determined by its first Chern class.
On the other hand,
$H^2(X,
 \zz)=H^2(M,\zz)=0$, since $M$ is a deformation retraction
of
$X$. Hence, the normal bundle of $X$ in $\cc^{n+1}$ is trivial. 
 Now, \rt{1} implies that $M$ admits a totally real immersion
in $\cc^n$. The proof
 of \nrc{2} is complete.

We notice that a real surface can be immersed in
$\rr^3$. Thus, \nrc{2} gives another proof that
all orientable surfaces admit totally real immersions
in $\cc^2$, which is  due to
Forstneri\v{c}~\cite{forstneric} when the surfaces are compact.
\nrc{2} is inconclusive when $M$
is a non-orientable compact surface, since $H^2(M,\zz)=\zz_2$.
It was proved  by Forstneri\v{c}~\cite{forstneric}  that a non-orientable
compact surface admits a totally real immersion in $\cc^2$  if and only if
its genus is even.
In view of \nrc{1}, we see that a totally real
immersion of a non-orientable compact surface in $\cc^3$
has a trivial complex normal bundle if and only if the genus of the surface is
even.  This also indicates that the hypothesis that $H^2
(M,\zz)=0$
is needed in \nrc{2}, although it is not a necessary condition for the
triviality of $T^cM$.

\setcounter{thm}{0}\setcounter{equation}{0}
\section{ Proof of \rt{1}}

In this section we shall  prove \rt{1} by 
using \rt{2}. We shall see that the proof of \rt{1} is eventually
related to homotopy groups of complex Stiefel manifolds.

Let us first consider the case $k\leq 4$. We need the following lemma.
\bl{fg}
Let $M$ and $N$ are two totally real and real analytic 
submanifolds of $\cc^n$. Let
 $\nu_M$ and $\nu_N$ be the
complex normal bundles of $M$ and $N$ respectively.
Then $M$ and $N$ are biholomorphically
equivalent if and only if there is an analytic diffeomorphism $f\colon M\to N$ such
that the pull-back $f^*\nu_N$ is topologically isomorphic to $\nu_M$.
\el
\begin{proof} If $M, N$ are biholomorphically equivalent by $F$, then 
 it is clear that $dF$ maps $ T^cM$ into $T^cN$
 and $df$ also induces an isomorphism from
$\nu_M$ to $\nu_N$.
We now assume
that there is an analytic diffeomorphism $f\colon M\to N$ such that $\nu_M$
is topologically isomorphic to $f^*\nu_N$. 
Let $f^c$ be an biholomorphic extension
of $f$ which sends a Stein neighborhood $M^c$ of $M$ onto $N^c=f^c(M^c)$. We may assume
that $M$ is a strong deformation retraction of $M^c$, 
i.e. there is a continuous mapping
$r\colon M^c\to M^c$ with $r|_M=\mbox{id}$ and a homotopy $r_t\colon M^c\to M$
with $r_0=r$ and
$r_1=\mbox{id}$.  This implies that
any complex vector bundle $V$ over $M^c$ is isomorphic to $r^*(V|_M)$.
Therefore, $\nu_{M^c}$ is topologically isomorphic to $f^*\nu_{N^c}$. By a theorem of Grauert~\cite{grauertfiber},
$\nu_{M^c}$ and $\nu_{N^c}$ are also holomorphically isomorphic.
Let $s$ be the zero section of $\nu_{M^c}$.
By a theorem of Docquier and Grauert (see~\cite{gunning}, p.~257), the zero-section mapping can be extends to
a biholomorphic mapping from a neighborhood of $M^c$ in $\cc^n$ to
a neighborhood of zero section of $\nu_{M^c}$.
Therefore, $f^c$ extends to biholomorphic
mapping defined in a neighborhood of $M$ in $\cc^n$. %mention FR
\end{proof}

To start the proof of \rt{1}, we first notice that Forstneri\v{c} and Rosay~\cite{fr}
proved that
any two totally real
immersions
of a smooth manifolds $M$ in $\cc^n$ are regularly homotopic through totally
real immersions,
if $\text{dim} M\leq 2n/3$. Consequently, their argument also
 showed that such two totally real and real analytic
embeddings are also biholomorphically equivalent if $M$ is compact. In the special case of spheres, we notice
that 
 the standard embedding of $S^k$ in $\cc^{k+1}$ 
has a trivial complex normal bundle, 
then \rl{fg} and \nrc{normal} also imply that any two totally real and real 
analytic  embeddings of $S^k$ in $\cc^n$ 
are  biholomorphically equivalent when $n\geq 3k/2$. In particular, all totally
real and real analytic embeddings of $S^k$ in $\cc^n$ are biholomorphically
equivalent when $k=1, 2$.
We also notice that  complex line bundles
on $S^k$ is classified by the group $\pi_{k-1}(U_1)$ (see~\cite{bottbk}).
Hence, complex line bundles on $S^k$ are trivial if $k\neq 2$. 
Now it is easy to see that
 all totally real and real analytic embeddings of $S^k$
in $\cc^n$
 are biholomorphically
equivalent if  $k\leq 4$.

We now consider the case of $k>4$. Assume first that for some $n>k$,
all totally real and real
 analytic embeddings of $S^k$ in $\cc^n$ are
biholomorphically equivalent.
This implies that all totally real and real analytic embeddings of
$S^k$ in $\cc^n$ have trivial
complex normal bundle.
Notice that when $n>k$,  a $C^1$-smooth totally real immersions of $S^k$ in $C^n
$
can be connected to a totally real and real analytic embeddings of $S^k$ in
$\cc^n$  by $C^1$ totally real immersions.  
Thus, \nrc{normal} implies that for any totally real immersion of $S^k$ in
$\cc^n$, its  complex normal bundle
 is trivial. Therefore, the commutative diagram \re{can} and
\rt{2} imply that 
\be{notonto}
i_*\colon \pi_k(U_k)\to\pi_k(U_{n,k})
\end{equation}
is an epimorphism, where $i_*$ is induced by the inclusion $U_k\subset U_{n,k}$.
Therefore, the proof of \rt{1} will be complete if we can show
that \re{onto} is not an epimorphism
for $n=n_k$, and also that
\be{onto}
i_*\colon\pi_k(U_{n_k,k})\to\pi
_k(U_{n,k})
\end{equation}
is an epimorphism for all $n>n_k$. Notice that
$\pi_k(U_{n,k})=0$ for $n\geq 3k/2$. This  also follows from the
fact that
all totally real immersions of $S^k$ in $\cc^n$ are regularly homotopic.
Hence, it suffices to show that \re{onto} holds for $n_k<n< 3k/2$.

It is well-known from the Bott periodicity theorem that $\pi_{2l}(U_n)=0$
 and
$\pi_{2l+1}(U_n)=\pi_{2l+1}(U_{l+1})=\zz$ for $n>l$.
We shall discuss in two cases.

\vspace{1.5ex}
\noindent
{\bf Case 1.} $k=2l\ (l\geq 3)$.
In this case,
 $\pi_{2l}(U_{2l})=0$ implies that
$n_k$ is the largest integer $n$ such
that $\pi_{2l}(U_{n,2l})$ is non-trivial. From [14, p.~127], one sees that
$$
\pi_{2l}(U_{3l-1,2l})=
\left\{\begin{array}{ll}
0, & \text{ if $l$ is even and $l\geq 2$}, \vspace{.75ex}\\  
\zz_2, & \text{ if $l$ is odd and $l\geq 3$}.
\end{array}\right.
$$
Hence, $n_{2l}=3l-1$ when $l$ is odd and $l\geq 3$. We now assume that
$l$ is even. Then, one has
 $$
\pi_{2l}(U_{3l-2,2l})=\left\{
\begin{array}{ll}
\zz_2, &l=4, \vspace{.75ex}\\  
\zz_{48/U(l+1,3)},& l\geq 6,
\end{array}\right. $$
where $U(\cdot, 3)$ is a James number which divides $24$
[14, p.~128]. In particular,
 $\pi_{2l}(U_{3l-2,2l})$ $\neq 0$, if $l$ even and
$l\geq 4$. This showed that if $l$ is even and $l\geq 4$, $n_{2l}=3l-2$
satisfies the property stated in \rt{1}.

\vspace{1.5ex}
\noindent {\bf Case 2.} $k=2l+1\ (l\geq 2)$. In this case,
one has
$\pi_{2l+1}(U_{2l+1})=\zz$. From~\cite{sigrist}, we find
\begin{gather}
\label{3l1} \pi_{2l+1}(U_{3l+1,2l+1})  = \zz,\\
\label{3l} \pi_{2l+1}(U_{3l,2l+1}) =\left\{ 
\begin{array}{ll}
\zz, &\text{ if $l$ is even and $l\geq 2$}, \vspace{.75ex}\\
\zz+\zz_2, &\text{ if $l$ is odd and $l\geq 3$}.
\end{array}\right.
\end{gather}
We now consider the homomorphism
\be{i}
i_*\colon \pi_{2l+1}(U_{3l,2l+1})\to\pi_{2l+1}(U_{3l+1,2l+1}),
\end{equation}
where $i_*$ is induced by the inclusion $U_{3l,2l+1}\subset U_{3l+1,2l+1}$.
We need the following.
\bl{i}
Let $i_*$ be defined by $\re{i}$. 
Then $i_*$  is an epimorphism if and only if $l$ is odd.
\el
Let us postpone the proof of \rl{i} for a while and finish our proof of
\rt{1}.
Assume first
 that $l$ is odd.
It is clear that there is no epimorphism from $\zz$ to $\zz+\zz_2$.
Hence, \re{notonto} is not onto for $k=2l+1$ and $n=3l$. On the other hand,
\rl{i} and the vanishing of $\pi_k(U_{n,k})$ for $n>3l+1$ imply  that \re{onto}
is onto for all $n>3l$. Therefore, $n_{2l}=3l$ satisfies the property stated
in \rt{1}.
Next, we assume that $l$ is even.
Then \rl{i} implies that
\re{notonto} is not onto for $n=3l+1$.
On the other hand, \re{onto} is onto for $n=3l+2$ and $k=2l+1$ because of
the vanishing of $\pi_{2l+1}(U_{3l+2,2l+1})$. Therefore,  we
conclude that
for even $l$, $n_{
2l+1}=3l+1$ satisfies the property stated in \rt{1}.

We now turn to the proof of \rl{i}. Let $
U_{3l}\to U_{3l,2l+1}$ be the standard
fibration with fiber $U_{l-1}$,
and $U_{3l+1}\to U_{3l+1,2l+1}$  the fibration with fiber $U_{l}$.
Then the inclusion
$U_{3l,2l+1}\subset U_{3l+1,2l+1}$ induces
 the following commutative
diagram of exact sequences
$$ \begin{array}{cccccc} %line one
 \pi_{2l+1}(U_{3l}) & \stackrel{j_*}{\longrightarrow}         & \pi_{2l+1}(U_{3l,2l+1})
          & \longrightarrow    &  \pi_{2l}(U_{l-1}) & \rightarrow
                                           \pi_{2l}(U_{3l}) \\ %line two
 \| &  &\downarrow i_* & &
\downarrow i_*^\prime & \| \\ %line three
 \pi_{2l+1}(U_{3l+1}) & \stackrel{j_*^\prime}{\longrightarrow}         & \pi_{2l+1}(U_{3l+1,
2l+1})    & \stackrel{\delta}{\longrightarrow}    &  \pi_{2l}(U_l) &
\rightarrow \pi_{2l}(U_{3l+1}). \\ %line four
\| & & & &  \downarrow p_*  & \| \\ %line five
 \zz & & & & \pi_{2l}(S^{2l-1}) & 0
\end{array} $$
It is clear that if $i_*$ is epimorphic, so is $i_*^\prime$.
However, one knows that
$\pi_{2l}(U_l)=\zz_{l!}$ (see \cite{bott}), and $\pi_{2l}(U_{l-1})=\mathbb
Z_{l!/2}$ for $l$ even (see~\cite{kervaire}). Thus $i_*^\prime $ is not epimorphic
when $l$ is even.  We now assume that $l$ is odd.  In this case, Kervaire
showed that $p_*=0$ (see~\cite{kervaire}, Lemma I.1).
Hence,
$i_*^\prime$ is an
epimorphism, i.e. $\delta i_*$ is an epimorphism. Using (\ref{3l1}) and
(\ref{3l}), we can write more explicitly the following diagram
$$\begin{array}{cccc}
\zz & \stackrel{j_*}{\longrightarrow}  & \zz+\zz_2 &  \\
\| &  &  \downarrow i_*  & \\
\zz &\stackrel{j_*^\prime}{\longrightarrow}
 & \zz & \stackrel{\delta}{\longrightarrow}
 \zz_{l!}\rightarrow 0.
\end{array}$$
Write $j_*(1)=g+\epsilon$, where $g\in \zz$ and $\epsilon
\in \zz_2$.
Clearly, $i_*(\epsilon)=0$. Hence, we get
$$
g\cdot i_*(e)=i_*j_*(1)=j_*^\prime(1)=\pm l!
$$
for $e=1\in \zz\subset\zz+\zz_2$. Thus, $i^*(e)$ divides $l!$.
On the other hand,
$\delta\circ i_*$ is an epimorphism. Hence, $i_*(e)\delta(1)$
must be a generator of $\zz_{l!}$.  Therefore $i_*(e)=\pm 1$, i.e.
$i_*$ is an epimorphism. The proof of \rl{i} is complete.

Next, we want to show that all totally real and real analytic embeddings
of $S^k$ in $\cc^n$ are unimodularly equivalent if $n>k$ and $k\leq 4$,
i.e. they are equivalent through a biholomorphic mapping $F$ satisfying
$F^*\Omega=\Omega$.
 To this end we shall
use the Cauchy-Kowalewski theorem to prove a slightly general
result.
\bp{unimodular}
Let $M, N$ be two totally real and real analytic submanifolds of $\cc^n$ which are biholomorphically equivalent. Assume
 that for some complexification $M^c$,
the holomorphic normal bundle of $M^c$ contains a holomorphic
subbundle of rank one. Then $M$ and $N$ are unimodularly
equivalent.
\ep
\begin{proof} Assume that $M$ and $N$ are equivalent by a biholomorphic mapping
$\varphi$
defined near $M$. It suffices to show that there is a biholomorphic mapping
$\psi$ defined near $M$
such that $\psi^*\Omega=\varphi^*\Omega$ and
$\psi(M)=M$. By shrinking $M^c$ if necessary, we may assume that $M^c\subset
\cc^n$ is a
Stein manifold. Now we have the decomposition
$\nu_{M^c}=\nu'\oplus L$, where
$L$ is a line bundle. Moreover,
we may
assume that $\nu$ is a subbundle
of $M^c\times \cc^n$. Thus, a
neighborhood $U$ of the zero section of $\nu_{M^c}$ is identified with
a neighborhood of $M^c$.
We now want to show that there is a holomorphic mapping
$$
\psi\colon (u,v)\mapsto (u,\lambda(u,v)v),\quad u\in\nu', \ v\in L
$$
such  that $\psi^*\Omega=\varphi^*\Omega$.
   Let $w=(w_1,\cdots, w_k)$ be local coordinates on $M^c$, and let
 $\xi=(\xi_1,\ldots,\xi_{n-k-1})$ and $t$ be local trivializations of
$\nu'$ and $L$ respectively.
 Then
in local coordinates,
 $\psi$ must be in the form $(w,\xi,t)\to (w,\xi, t'(w,\xi,t))$ with
 $t'|_{t=0}=0$.
 Put
 $$
\varphi^*\Omega=f(w,\xi,t)dw\wedge d\xi\wedge dt,\quad
 \Omega=a(w,\xi,t)dw\wedge d\xi\wedge dt,
 $$
 where $dw=dw_1\wedge\ldots\wedge w_k$
, $d\xi=d\xi_1\wedge\ldots\wedge d\xi_{n-k-1}$.
 Then $\psi^*\Omega=\varphi^*\Omega$ is equivalent to the equation
 $a(w,\xi,t'){\partial t'}/{\partial t}=f(w,\xi,t)$. By the
Cauchy-Kowalewski
 theorem, we know that for small $|t|$ the solution $t'$ exists
uniquely. This means
 that the required mapping $\psi$ is uniquely determined
 in local coordinates. Therefore, there is a holomorphic mapping defined
 in a neighborhood of $M^c$ such that
 $\psi^*\Omega=\varphi^*\Omega$. Since the restriction of $\psi$ to
$M^c$ is
 the identity mapping, it is easy to see that $\psi$ is one-to-one in some 
neighborhood
 of $M^c$.
\end{proof}

%We remark that when $M$ is a sphere,
% the assumption in~\rp{unimodular} about
%the existence of line subbundle 
%is equivalent to the complex normal bundle of the sphere
% has a trivial line subbundle.
Notice that all totally real and real analytic embeddings of a compact
surface $M$ in $\cc^n\ (n\geq 3)$ are biholomorphically equivalent.
 Also, $M$ admits a totally real and real analytic embedding in $\cc^3$. Hence, the complex normal bundle of a totally real embeddings of
$M$ in $\cc^n\ (n\geq 3)$ is the direct sum of a line bundle and a trivial
bundle. From \rp{unimodular},
we have the following.
\bc{uni}
All totally real and real analytic embeddings of a
compact surface in
$\cc^n\ (n\geq 3)$ are unimodularly equivalent.
\ec

 \bibliographystyle{plain}

\end{document}